\documentclass[11pt]{amsart}

\usepackage[a4paper,margin=1.08in]{geometry}
\usepackage{amsmath,amssymb,mathtools}
\usepackage{microtype}
\usepackage{needspace}
\usepackage{enumitem}
\usepackage{booktabs}
\usepackage{xcolor}
\usepackage[colorlinks=true,linkcolor=blue!55!black,citecolor=blue!55!black,urlcolor=blue!55!black]{hyperref}

\newtheorem{theorem}{Theorem}[section]
\newtheorem{lemma}[theorem]{Lemma}
\newtheorem{corollary}[theorem]{Corollary}
\newtheorem{proposition}[theorem]{Proposition}
\theoremstyle{definition}

\theoremstyle{remark}
\newtheorem{remark}[theorem]{Remark}

\newcommand{\F}{\mathbb F}
\newcommand{\calF}{\mathcal F}

\newcommand{\ex}{\operatorname{ex}}

\newcommand{\E}{\mathbb E}
\newcommand{\Prob}{\mathbb P}
\newcommand{\betafive}{\beta_5}

\title[An improved bound for the Tur\'an densities of daisies]{An improved asymptotic lower bound for the Tur\'an densities of daisies}
\author{Ziyi Feng}
\address[Ziyi Feng]{School of Mathematical Sciences, Anhui University,
111 Jiulong Road, Hefei, Anhui 230601, China}
\email{A125211046@stu.ahu.edu.cn}
\date{November 25, 2026}

\begin{document}

\begin{abstract}
Ellis, Ivan and Leader constructed daisy-free hypergraphs by retaining
sets whose vector labels form a basis over $\F_2$.  We extend their
construction by labelling vertices with nonzero vectors of $\F_q^r$
and retaining an $r$-set when its determinant belongs to a prescribed
sign-symmetric set $A\subseteq\F_q^\times$.  The binary basis
construction is the special case $q=2$, $A=\{1\}$.  Over $\F_5$,
accepting the values $\{\pm1\}$ improves the limiting lower bound for
the Tur\'an density of daisies from $\beta_2=0.288788\ldots$ to
$\beta_5/2=0.380166\ldots$, where
$\beta_q=\prod_{k\ge1}(1-q^{-k})$.

We prove that a nonempty sign-symmetric set $A$ guarantees
daisy-freeness for every labelling and every dimension if and only if
$(AA+AA)\cap AA=\varnothing$, where $AA=\{ab:a,b\in A\}$.
Under independent uniform sampling, each fixed pair $(q,A)$ gives
limiting expected density $\beta_q|A|/(q-1)$.  We show that its maximum
over all such rules is $\beta_5/2$, attained precisely when $q=5$
and $A=\{\pm a\}$ for some $a\in\F_5^\times$.
Uniform convergence shows that dimension-dependent choices do not
improve this constant.

Applications give dense induced $C_6$-free
subgraphs of Boolean edge-layers and vertex sets of the hypercube
meeting every four-dimensional subcube in at most five vertices.
\end{abstract}

\maketitle

\section{Introduction}

For $r\ge2$, the $r$-\emph{daisy} $D_r$ is the $r$-uniform hypergraph
\[
 D_r=\{S\cup X:X\in\tbinom T2\},
\]
where $|S|=r-2$, $|T|=4$, and $S\cap T=\varnothing$.
The set $S$ is called the \emph{stem}.  For an $r$-uniform hypergraph $H$,
let $\ex(n,H)$ denote the maximum number of edges in an $H$-free
$r$-uniform hypergraph on $n$ vertices.  We denote its Tur\'an density by
\[
 \pi(H)=\lim_{n\to\infty}\frac{\ex(n,H)}{\binom nr}.
\]
For the
sequence of daisies $D_r$, we subsequently let $r\to\infty$; thus the
quantity of interest is the  limit
\[
 \lim_{r\to\infty}\pi(D_r)
 =\lim_{r\to\infty}\left(
     \lim_{n\to\infty}\frac{\ex(n,D_r)}{\binom nr}
   \right).
\]
Bollob\'as, Leader and Malvenuto \cite{BLM},
and independently Bukh (as reported in \cite{BLM}), conjectured that
$\pi(D_r)\to0$ as $r\to\infty$.
Ellis, Ivan and Leader \cite[Theorem~4]{EIL} disproved this conjecture using a basis
construction over $\F_2$ and obtained
\begin{equation}\label{eq:EIL}
 \lim_{r\to\infty}\pi(D_r)\ge
 \prod_{k\ge1}(1-2^{-k})=0.288788\ldots.
\end{equation}

To extend their construction, we begin with the following observation.
For an ordered list of column vectors $x_1,\dots,x_r\in\F_q^r$,
we use the determinant notation
\[
 [x_1,\dots,x_r]
 :=\det\begin{pmatrix}x_1&\cdots&x_r\end{pmatrix}.
\]
For vertex labels $v_1,\dots,v_n$ and an $r$-set
$S=\{i_1<\cdots<i_r\}$, we denote their determinant in increasing
index order by
\[
 \Delta(S):=[v_{i_1},\dots,v_{i_r}].
\]
Over $\F_2$, since $1$ is the only nonzero field element, the labels of $S$
form a basis if and only if $\Delta(S)=1$.  Thus the
Ellis--Ivan--Leader construction can be viewed as accepting precisely
those $r$-sets whose determinant equals $1$.

Over any finite field $\F_q$, fix the ordered stem columns
$w_1,\dots,w_{r-2}$.  Throughout the paper, when these common columns
are fixed, we use the shorthand
\[
 [xy]:=[w_1,\dots,w_{r-2},x,y].
\]
For four outer vectors $a,b,c,d$, these brackets satisfy the
Grassmann--Pl\"ucker identity
\[
 [ab][cd]-[ac][bd]+[ad][bc]=0.
\]
The identity holds for any assignment of vectors to this common-stem
configuration; a daisy in the selected family would additionally
require all six determinants to be accepted.  We derive the identity
and explain the role of the common stem in Section
\ref{sec:construction}.  Over $\F_2$, six accepted values would all
equal $1$, contradicting the identity.

This observation suggests choosing a different field and a set of
accepted determinant values that still excludes such a configuration.
We work over $\F_5$ and retain an $r$-set when its determinant belongs
to $\{\pm1\}$.  Accepting both signs makes this condition independent
of the column order.  If all six petals of a daisy were retained,
each of the three signed terms in the identity would lie in
$\{\pm1\}$.  Their sum would be one of $\pm1,\pm3$, all nonzero
in $\F_5$, giving a contradiction.  This choice therefore remains
daisy-free and, under independent uniform sampling, yields the
improved limiting density quantified below.
We also ask whether this choice gives the largest asymptotic expected
density under independent uniform sampling, allowing any finite field
and any symmetric set of accepted determinant values that guarantees
daisy-freeness for every labelling.  We formulate this optimization
problem precisely after Theorem~\ref{thm:daisy-main}.

For each prime power $q$, we use the notation
\begin{equation}\label{eq:beta5}
 \beta_q:=\prod_{k\ge1}(1-q^{-k}),
 \qquad
 \betafive=0.760332795871232\ldots.
\end{equation}
The constant $\beta_q$ has a natural probabilistic meaning.  If we
choose $v_1,\dots,v_r$ independently and uniformly from
$\F_q^r\setminus\{0\}$, then
\[
 \Prob(v_1,\dots,v_r\text{ form a basis})
 =\frac{\prod_{k=1}^{r}(1-q^{-k})}{(1-q^{-r})^r}
 \longrightarrow\beta_q
 \qquad(r\to\infty),
\]
with $q$ fixed.  Thus $\beta_q$ measures the limiting proportion of
ordered samples that form a basis.  Conditional on this event, the
determinant is uniformly distributed over the $q-1$ nonzero field
elements.  In particular, over $\F_5$, accepting only the values
$\{\pm1\}$ retains half of these samples, which explains the limiting
edge density $\beta_5/2$ in the following theorem.  We give the
counting argument in Section \ref{sec:construction}.
Our first result records both the construction and the lower bound it gives.

\begin{theorem}[The $\F_5$ construction]\label{thm:daisy-main}
Let $2\le r\le n$, and assign nonzero vectors $v_1,\dots,v_n\in\F_5^r$
to the elements of $[n]$.  The family
\begin{equation}\label{eq:Fdef}
 \calF=\left\{S\in\binom{[n]}r:
                 \Delta(S)\in\{\pm1\}\right\}
\end{equation}
is $D_r$-free for every choice of labels.  If $v_1,\dots,v_n$ are chosen
independently and uniformly from $\F_5^r\setminus\{0\}$, then
\begin{equation}\label{eq:F-density}
 \begin{aligned}
 \E|\calF|
 &=\frac12\,
   \frac{\prod_{k=1}^r(1-5^{-k})}{(1-5^{-r})^r}\binom nr\\
 &>\frac{\betafive}{2}\binom nr.
 \end{aligned}
\end{equation}
In particular,
\[
 \lim_{r\to\infty}\pi(D_r)\ge
 \frac{\betafive}{2}=0.380166397935616\ldots.
\]
\end{theorem}

Thus the construction improves the lower bound of Ellis, Ivan and
Leader \cite[Theorem~4]{EIL}, stated in equation \eqref{eq:EIL}, from
\[
 0.288788\ldots\quad\text{to}\quad0.3801663979\ldots.
\]

We next ask whether the $\F_5$ construction is optimal within the class
of constructions obtained by prescribing determinant values over a
finite field.  To formulate the question precisely, we specify the
random construction, the requirement of daisy-freeness, and the density
to be optimized separately.

\smallskip
\noindent\textbf{The random construction.}
Let $q$ be any prime power and let $A\subseteq\F_q^\times$ be any
nonempty \emph{sign-symmetric} set, meaning $A=-A$.  For integers
$n\ge r\ge2$, we choose vectors $v_1,\dots,v_n$ independently and
uniformly from $\F_q^r\setminus\{0\}$.  We denote the resulting
random family by
\begin{equation}\label{eq:general-filter}
 \calF_A(v_1,\dots,v_n)
 =\left\{S\in\binom{[n]}r:\Delta(S)\in A\right\}.
\end{equation}
Here $n$ is the number of vertices and $r$ is both the edge size and
the dimension of the label space.  Repeated vectors are allowed, so
$n$ is not bounded by the number $q^r-1$ of possible labels.
When $(q,A)$ is fixed and $r$ varies, the coefficient field remains
$\F_q$: the determinant of an $r\times r$ matrix still takes values
in $\F_q$, so the same set $A$ defines an acceptance rule in every
dimension.  For an unordered $r$-set, permuting
the corresponding columns changes the determinant by at most a sign;
the condition $A=-A$ therefore makes membership independent of the
chosen ordering.  At this stage neither the field nor the set of
accepted determinant values is restricted beyond these assumptions.

\smallskip
\noindent\textbf{The daisy-free requirement.}
  For a fixed field $\F_q$, we call
$A$ \emph{nice} if
\[
 \begin{gathered}
 \calF_A(v_1,\dots,v_n)\text{ is }D_r\text{-free}\\
 \text{for every }n\ge r\ge2\text{ and every }
 (v_1,\dots,v_n)\in(\F_q^r\setminus\{0\})^n.
 \end{gathered}
\]
Equivalently, for every $n\ge r\ge2$, the random family above is
$D_r$-free with probability one: each possible assignment of vectors
has positive probability. We denote the class of nice pairs by
\[
 \mathcal U:=\left\{(q,A):
 \begin{array}{l}
 q\text{ is a prime power},\\
 \varnothing\ne A\subseteq\F_q^\times,\quad A=-A,\\
 A\text{ is nice}
 \end{array}\right\}.
\]
This class does not depend on $r$ or $n$.

\smallskip
\noindent\textbf{The optimization problem.}
For $r\ge2$ and $(q,A)\in\mathcal U$, we choose $v_1,\dots,v_r$
independently and uniformly from $\F_q^r\setminus\{0\}$ and denote
the probability of an accepted determinant by
\[
 \begin{aligned}
 p_r(q,A)
 &:=\Prob\bigl([v_1,\dots,v_r]\in A\bigr)\\
 &=\frac{|A|}{q-1}\,
   \frac{\prod_{k=1}^{r}(1-q^{-k})}{(1-q^{-r})^r}.
 \end{aligned}
\]
The product quotient is the probability that the labels form a basis;
conditional on this event, their determinant is uniform on
$\F_q^\times$ and is accepted with probability $|A|/(q-1)$.
We give the counting argument in Section \ref{sec:construction}.

For $n\ge r$, linearity of expectation gives
\[
 \E|\calF_A(v_1,\dots,v_n)|=\binom nr\,p_r(q,A),
\]
so the expected edge density is independent of $n$.  Every outcome
is $D_r$-free, and some outcome has at least the expected number of
edges.  Taking $n\to\infty$ with $r,q,A$ fixed therefore gives
\begin{equation}\label{eq:filter-turan-bound}
 \pi(D_r)\ge p_r(q,A).
\end{equation}

Each nice pair $(q,A)$ specifies a rule that gives a construction
in every dimension $r$.  We keep this pair fixed and let
$r\to\infty$, after the vertex-number limit above.  We denote the
limiting expected edge density by
\begin{equation}\label{eq:d-objective}
 d(q,A):=\lim_{r\to\infty}p_r(q,A)
 =\beta_q\frac{|A|}{q-1}.
\end{equation}
Subsection \ref{sec:densities} establishes this limit.  Each nice pair
therefore supplies a lower bound for the same limiting Tur\'an density:
\[
 \lim_{r\to\infty}\pi(D_r)\ge d(q,A).
\]
Our main question is to characterize $\mathcal U$ and determine
which fixed rules give the best such bound.  More precisely, we seek
\begin{equation}\label{eq:optimization-problem}
 \sup_{(q,A)\in\mathcal U}d(q,A)
 =\sup_{(q,A)\in\mathcal U}
   \left(\lim_{r\to\infty}p_r(q,A)\right),
\end{equation}
and all pairs attaining it.  The pair $(q,A)$ remains fixed in each
limit.  We optimize expected densities under independent uniform
sampling; maximizing over individual assignments of labels is a
different problem.

For comparison, Alon \cite[Corollary~2.4]{Alon} considered the binary
basis construction of Ellis, Ivan and Leader with labels drawn
independently from an arbitrary common distribution on
$\F_2^r\setminus\{0\}$.  Keeping the field and basis rule fixed, he
proved that the uniform distribution maximizes the probability of
obtaining a basis for each $r$.  Our optimization varies different
parameters: the distribution is uniform on all nonzero vectors, while
the field and accepted determinant set may vary.

To study \eqref{eq:optimization-problem}, we express the obstruction
to a daisy in the language of additive combinatorics.  The
Grassmann--Pl\"ucker relation governs additive relations among products
of accepted determinant values, leading to the product set
\[
 B=AA:=\{ab:a,b\in A\}.
\]
The following theorem characterizes rules that exclude daisies for
every labelling; Sections \ref{sec:construction} and \ref{sec:criterion}
derive the determinant relation and prove the criterion.

\begin{theorem}[The daisy-free criterion]\label{thm:criterion}
A nonempty sign-symmetric set $A\subseteq\F_q^\times$ is nice if and only if
\begin{equation}\label{eq:criterion}
 0\notin AA-AA+AA.
\end{equation}
Equivalently, for $B=AA$,
\begin{equation}\label{eq:sumfree}
 (B+B)\cap B=\varnothing.
\end{equation}
\end{theorem}

The two conditions have distinct roles.  The symmetry $A=-A$ ensures
that the edge definition does not depend on the ordering of its vector
labels.  In contrast, the additive condition $(AA+AA)\cap AA=\varnothing$
ensures daisy-freeness.

Having identified the nice pairs, we solve
\eqref{eq:optimization-problem} and determine all maximizing pairs.

\begin{theorem}[Asymptotic optimality for independent uniform samples]\label{thm:finite-field-optimal}
Among the fixed pairs $(q,A)\in\mathcal U$,
\[
 \max_{(q,A)\in\mathcal U}d(q,A)=\frac{\beta_5}{2}.
\]
The maximum is attained precisely when
\[
 q=5,\qquad A=\{a,-a\}
\]
for some $a\in\F_5^\times$.
\end{theorem}

The theorem identifies the rules with the best limiting expected
density in this framework.  It does not assert that the $\F_5$ filter
maximizes $p_r(q,A)$ for every fixed $r$.

\begin{remark}[Rules depending on the dimension]\label{rem:varying-parameters}
One can also choose a different nice pair for each $r$.  If we denote
the best expected density in dimension $r$ by
\begin{equation}\label{eq:best-uniform-bound}
 M_r:=\sup_{(q,A)\in\mathcal U}p_r(q,A),
\end{equation}
then
\[
 \frac{\beta_5}{2}\le M_r
 \le\frac{\beta_5}{2}+(r+1)2^{-r},
 \qquad
 \lim_{r\to\infty}M_r=\frac{\beta_5}{2}.
\]
Indeed, Lemma \ref{lem:uniform-density} gives
$0\le p_r(q,A)-d(q,A)\le(r+1)2^{-r}$ uniformly over $\mathcal U$.
Taking suprema and applying Theorem \ref{thm:finite-field-optimal}
yields the bounds.  Thus allowing the parameters to vary with $r$
does not improve the asymptotic constant obtained by a fixed rule.
\end{remark}

We next turn to applications in the hypercube.  Let $Q_n$ denote
the graph on $\{0,1\}^n$ in which two vertices are adjacent when they
differ in exactly one coordinate.  Identifying each vertex with a
subset of $[n]$, we denote by $L_r(n)$ the incidence graph between
$\binom{[n]}{r-1}$ and $\binom{[n]}r$, with adjacency given by inclusion.
This is the $r$th Boolean edge-layer of $Q_n$.  For a graph $G$, we
denote its number of edges by $e(G)$, and $C_k$ denotes a cycle of
length $k$.

In their proof that $C_{10}$ has positive
Tur\'an density in the hypercube, Grebennikov and Marciano \cite{GMv4}
isolated the auxiliary problem of finding a dense induced $C_6$-free subgraph
inside a single Boolean edge-layer $L_r(n)$.  Applying our determinant-value
obstruction on its two vertex levels gives the following improvement.

\begin{corollary}\label{cor:layer-main}
For every $r,n$ with $2\le r\le n$, there exists an induced $C_6$-free subgraph $G_r\subseteq L_r(n)$ satisfying
\[
 \begin{aligned}
 e(G_r)&>\frac{\betafive}{5}\,e(L_r(n))\\
 &=0.152066559174246\ldots\,e(L_r(n)).
 \end{aligned}
\]
\end{corollary}

Grebennikov and Marciano \cite[Theorem~1.3]{GMv4} proved the corresponding binary-basis bound
\[
 \begin{aligned}
 e(G_r)&>\frac12\prod_{k=1}^{\infty}(1-2^{-k})\,e(L_r(n))\\
 &=0.1443940475\ldots\,e(L_r(n)).
 \end{aligned}
\]
Thus Corollary \ref{cor:layer-main} gives a slight improvement on their result.

Johnson and Talbot \cite[Question~13]{JT} asked whether, for every
$\delta>0$ and every fixed integer $d\ge2$, all sufficiently large $n$ have
the following property: every family $\mathcal A\subseteq\{0,1\}^n$ of
density at least $\delta$ contains at least
$\binom d{\lfloor d/2\rfloor}$ points in some $d$-dimensional subcube.
Ellis, Ivan and Leader \cite[Theorem~8]{EIL} answered this question
negatively for every $d\ge4$.  In the first such case, $d=4$, the proposed
threshold is $\binom42=6$.

Thus, already for $d=4$, one is led to the quantitative problem of maximizing
the asymptotic density of a family that meets every $4$-dimensional subcube
in at most five points.  The binary daisy construction of Ellis, Ivan and
Leader gives the lower bound
\[
 \frac13\prod_{k\ge1}(1-2^{-k})=0.0962626983\ldots,
\]
while the elementary local-averaging upper bound is $5/16$.
Replacing the binary daisy family by Theorem \ref{thm:daisy-main} improves
the lower endpoint of this open density problem as follows.

\begin{corollary}\label{cor:q4}
For $n\to\infty$ there exists $\mathcal A\subseteq\{0,1\}^n$ such that every $4$-dimensional subcube contains at most five points of $\mathcal A$, while
\[
\begin{aligned}
\frac{|\mathcal A|}{2^n}&\ge \frac{\betafive}{6}-o(1)\\
&=0.126722132645205\ldots-o(1).
\end{aligned}
\]
\end{corollary}

Thus Corollary \ref{cor:q4} raises the previously known lower bound from
$0.0962626983\ldots$ to $0.1267221326\ldots$; the upper bound $5/16$
is unchanged; closing the gap remains open.

Further substitutions into the arguments of Ellis, Ivan and Leader,
and into the $C_{10}$ construction of Grebennikov and Marciano, are
discussed in Section \ref{sec:conclusion}.  We do not state these as
additional main results: the Ellis--Ivan--Leader applications reuse their
existing arguments, while the resulting $C_{10}$ bound is weaker than
Peji\'c's bound \cite{Pejic}.

\smallskip
\noindent\textbf{Organization of the paper.}\par\nobreak
Section \ref{sec:construction} proves that the $\F_5$ family is daisy-free
and computes its density, completing the proof of
Theorem \ref{thm:daisy-main}.  Section \ref{sec:criterion} proves the exact
criterion in Theorem \ref{thm:criterion}.  Section \ref{sec:optimality}
proves Theorem \ref{thm:finite-field-optimal} by computing the limiting
expected density, applying a character estimate over arbitrary finite
fields, and completing the small-field comparisons.  Section \ref{sec:applications} proves
the two application corollaries and establishes optimality for the
two-level construction within the same framework of nice filters.
Section \ref{sec:conclusion} discusses
the scope of the optimality result and possible extensions.
Subsection \ref{sec:densities} includes the uniform convergence estimate
used in Remark \ref{rem:varying-parameters} to handle rules that
vary with the dimension.

\smallskip
\noindent\textbf{Notation.}
Throughout the paper, $[n]$ denotes $\{1,\dots,n\}$, and
$\binom Xr$ denotes the family of all $r$-element subsets of a finite
set $X$.  We denote the field of order $q$ by $\F_q$ and its
multiplicative group of nonzero elements by $\F_q^\times$, where $q$
is a prime power.  The symbols $\Prob$ and $\E$ denote probability
and expectation, respectively.

For subsets $X,Y$ of a field, we use the standard sumset,
difference-set, and product-set notation
\[
 X+Y=\{x+y:x\in X,\ y\in Y\},\qquad
 X-Y=\{x-y:x\in X,\ y\in Y\},
\]
\[
 XY=\{xy:x\in X,\ y\in Y\}.
\]

\section{The \texorpdfstring{$\F_5$}{F5} family is daisy-free}\label{sec:construction}

\begin{proof}[Proof of Theorem \ref{thm:daisy-main}]
Membership in \eqref{eq:Fdef} does not depend on the column order,
since both $1$ and $-1$ are accepted.

\smallskip
\noindent\emph{The identity supplied by a common stem.}
For use here and in the general criterion, we work over an arbitrary finite
field.  We fix $r-2$ vectors $w_1,\dots,w_{r-2}$ and four further vectors
$a,b,c,d$.  Using the common-stem bracket notation introduced above,
we prove the identity
\begin{equation}\label{eq:plucker}
 [ab][cd]-[ac][bd]+[ad][bc]=0.
\end{equation}
First suppose that the common vectors $w_1,\dots,w_{r-2}$ are linearly
dependent.  Every determinant $[xy]$ then has linearly dependent columns
and is zero, so the identity holds.  We may therefore assume that these
vectors are linearly independent. Let
\[
 U=\operatorname{span}(w_1,\dots,w_{r-2}),
\]
with $U=\{0\}$ when $r=2$.  If $a\in U$, then $[ab]=[ac]=[ad]=0$,
so again every term in \eqref{eq:plucker} is zero.

It remains to consider $a\notin U$.  We extend the linearly independent
list $w_1,\dots,w_{r-2},a$ to a basis by adding a vector $e$, and set
\[
 \lambda=[w_1,\dots,w_{r-2},a,e]\ne0.
\]
We write the remaining vectors in this basis as
\[
 \begin{aligned}
 b&=u_b+b_1a+b_2e,\\
 c&=u_c+c_1a+c_2e,\\
 d&=u_d+d_1a+d_2e,
 \end{aligned}
 \qquad u_b,u_c,u_d\in U.
\]
The components in $U$ do not contribute to any bracket: subtracting
a linear combination of the first $r-2$ columns from either of the last
two columns leaves the determinant unchanged.  Thus only the two
coordinates along $a$ and $e$ matter.  More explicitly, if
$x=u_x+x_1a+x_2e$ and $y=u_y+y_1a+y_2e$, then multilinearity and
alternation of the determinant give
\[
 \begin{aligned}
 [xy]
 &=x_1y_2[w_1,\dots,w_{r-2},a,e]
   +x_2y_1[w_1,\dots,w_{r-2},e,a]\\
 &=\lambda(x_1y_2-x_2y_1).
 \end{aligned}
\]
The omitted terms have two equal final columns and hence vanish.
Since $a$ has coordinates $(1,0)$, the six brackets are therefore
\[
 \begin{aligned}
 [ab]&=\lambda b_2,& [cd]&=\lambda(c_1d_2-c_2d_1),\\
 [ac]&=\lambda c_2,& [bd]&=\lambda(b_1d_2-b_2d_1),\\
 [ad]&=\lambda d_2,& [bc]&=\lambda(b_1c_2-b_2c_1).
 \end{aligned}
\]
Substitution into the left-hand side of \eqref{eq:plucker}, followed by
division by $\lambda^2$, gives the following pairwise cancellation:
\begin{align*}
 &b_2(c_1d_2-c_2d_1)-c_2(b_1d_2-b_2d_1)
   +d_2(b_1c_2-b_2c_1)\\
 &\quad=(b_2c_1d_2-d_2b_2c_1)
       +(c_2b_2d_1-b_2c_2d_1)\\
 &\qquad\quad +(d_2b_1c_2-c_2b_1d_2)=0.
\end{align*}

\smallskip
\noindent\emph{Excluding a daisy.}
We now return to $\F_5$.  Suppose that $\calF$ contains a daisy with stem $I$
and four outer vertices $a,b,c,d$.  We order the common $I$-columns first.
Since each of the six petals belongs to $\calF$, sign symmetry gives
\[
 [ab],[ac],[ad],[bc],[bd],[cd]\in\{\pm1\}.
\]
Each of the three signed terms on the left of \eqref{eq:plucker} is
therefore in $\{\pm1\}$.  Their sum can only be one of $\pm1,\pm3$,
none of which is zero in $\F_5$.  This contradicts \eqref{eq:plucker},
so the family is daisy-free for every choice of labels.

\smallskip
\noindent\emph{Computing the expected density.}
We carry out the count over a general finite field $\F_q$ and then
specialize to $q=5$.  Let $\mathrm{GL}_r(q)$ denote the group of
invertible $r\times r$ matrices over $\F_q$, and let
$\mathrm{SL}_r(q)$ denote its subgroup of matrices with determinant $1$.
If the first $j$ columns are linearly independent, their span contains
$q^j$ vectors, so there are $q^r-q^j$ choices for the next column.
Thus
\[
 |\mathrm{GL}_r(q)|=\prod_{j=0}^{r-1}(q^r-q^j).
\]
There are $(q^r-1)^r$ equally likely ordered samples of $r$
independent uniform nonzero vectors.  The probability that such a
sample forms a basis is therefore
\[
 \Prob(v_1,\dots,v_r\text{ form a basis})
 =\frac{|\mathrm{GL}_r(q)|}{(q^r-1)^r}
 =\frac{\prod_{k=1}^r(1-q^{-k})}{(1-q^{-r})^r}.
\]
Conditional on this event, all ordered bases are equally likely.
For any $s,t\in\F_q^\times$, multiplying the first column by $t/s$
is a bijection between bases with determinant $s$ and bases with
determinant $t$.  Hence the determinant is uniform on
$\F_q^\times$, and each nonzero value occurs exactly
$|\mathrm{SL}_r(q)|$ times, where
\begin{equation}\label{eq:SL}
 |\mathrm{SL}_r(q)|
 =\frac{|\mathrm{GL}_r(q)|}{q-1}
 =\frac{q^{r^2}}{q-1}\prod_{k=1}^r(1-q^{-k}).
\end{equation}
For any nonempty $A\subseteq\F_q^\times$, it follows that
\begin{equation}\label{eq:random-det-prob}
 p_r(q,A)=\Prob\bigl([v_1,\dots,v_r]\in A\bigr)
 =\frac{|A|\,|\mathrm{SL}_r(q)|}{(q^r-1)^r}.
\end{equation}
Equivalently, the probability of retaining an $r$-set is the
probability of forming a basis multiplied by the conditional
probability that its determinant belongs to $A$:
\begin{equation}\label{eq:finite-daisy}
 p_r(q,A)
 =\underbrace{\frac{\prod_{k=1}^r(1-q^{-k})}{(1-q^{-r})^r}}
     _{\text{probability of a basis}}
 \times
 \underbrace{\frac{|A|}{q-1}}
     _{\substack{\text{fraction of accepted}\\
                  \text{nonzero determinant values}}}.
\end{equation}
When $A=-A$, every $r$-set is retained with this same probability,
so linearity of expectation gives
\[
 \E|\calF_A(v_1,\dots,v_n)|=\binom nr\,p_r(q,A).
\]
Thus $p_r(q,A)$ is also the expected edge density of the whole
family, independently of $n$.  The retention events for different
$r$-sets need not be independent.

For a fixed $S\in\binom{[n]}r$, we take $q=5$ and $A=\{\pm1\}$.
Equations \eqref{eq:SL} and \eqref{eq:random-det-prob} yield
\[
 \Prob(S\in\calF)
 =\frac{2|\mathrm{SL}_r(5)|}{(5^r-1)^r}
 =\frac12\,\frac{\prod_{k=1}^r(1-5^{-k})}{(1-5^{-r})^r}.
\]
Summing over all $S$ proves \eqref{eq:F-density}.  The finite product
exceeds $\betafive$ and the denominator is less than one, giving the
strict inequality there.  Since every outcome is daisy-free, some
outcome has at least the expected number of edges.  Letting $n\to\infty$
and then $r\to\infty$ gives the asserted lower bound for
$\lim_{r\to\infty}\pi(D_r)$.
\end{proof}

Repeated vector labels are allowed throughout this construction.  Any
repetition within a candidate edge makes its determinant zero and
therefore excludes that edge.

\section{Universal daisy-freeness and the Pl\"ucker obstruction}\label{sec:criterion}

We now explain why the general determinant-value construction leads to
an additive-combinatorial condition on the product set $AA$.  We consider
a  daisy and denote the six determinants associated with its
petals by $x_{ij}$, where $1\le i<j\le4$.  If all six petals are retained,
then every $x_{ij}$ belongs to $A$.  The Grassmann--Pl\"ucker relation is
\[
 x_{12}x_{34}-x_{13}x_{24}+x_{14}x_{23}=0.
\]
Each of its three terms, before the displayed signs are applied, is a
product of two elements of $A$ and hence belongs to $B=AA$.  Thus a
daisy would produce elements $b_1,b_2,b_3\in B$ satisfying
\[
 b_1-b_2+b_3=0,
\]
or equivalently $b_1+b_3=b_2$.  This identifies the relevant
obstruction: the product set $B$ must contain no solution to
$b_1+b_3=b_2$, which is precisely the sum-free condition
\[
 (B+B)\cap B=\varnothing,
\]
or, in terms of $A$, $0\notin AA-AA+AA$.  The sufficiency of this
condition follows immediately from the Pl\"ucker relation.  To prove
necessity, we show that every solution of the six-value relation in
$A$ can be realized as the six minors of four actual vectors.

\begin{proof}[Proof of Theorem \ref{thm:criterion}]
Suppose a daisy occurs, with ordered stem columns
$w_1,\dots,w_{r-2}$ and four remaining columns $a,b,c,d$.
All six common-stem brackets lie in $A$, while \eqref{eq:plucker} gives
\[
 [ab][cd]-[ac][bd]+[ad][bc]=0.
\]
Thus \eqref{eq:criterion} rules out every daisy.

Conversely, suppose that \eqref{eq:criterion} fails.  Then there are six
elements $x_{12},x_{13},x_{14},x_{23},x_{24},x_{34}\in A$ satisfying
\begin{equation}\label{eq:plucker-six}
 x_{12}x_{34}-x_{13}x_{24}+x_{14}x_{23}=0.
\end{equation}
In $\F_q^2$, we choose the vectors
\[
 u_1=\binom10,\qquad
 u_2=\binom0{x_{12}},\qquad
 u_3=\binom{-x_{23}/x_{12}}{x_{13}},\qquad
 u_4=\binom{-x_{24}/x_{12}}{x_{14}}.
\]
Their six pairwise determinants are exactly the six values $x_{ij}$; the
last equality, $[u_3,u_4]=x_{34}$, is \eqref{eq:plucker-six}.
For any $r\ge2$, we take the standard basis vectors
$e_1,\dots,e_{r-2}$ as the ordered stem columns.  If
$u_i=(\alpha_i,\beta_i)^{\mathsf T}$, we assign to its outer vertex
\[
 \widetilde u_i=\alpha_i e_{r-1}+\beta_i e_r.
\]
All labels are nonzero, and for every $i<j$,
\[
 [e_1,\dots,e_{r-2},\widetilde u_i,\widetilde u_j]=x_{ij}.
\]
Thus all six petals are retained, so $A$ is not nice.

Finally, $A=-A$ implies $B=-B$.  Therefore $0\in B+B+B$ if and only if
$(B+B)\cap B\ne\varnothing$, which proves the equivalence of
\eqref{eq:criterion} and \eqref{eq:sumfree}.
\end{proof}

\section{Why \texorpdfstring{$\F_5$}{F5} is asymptotically optimal}\label{sec:optimality}

We use the density formula from Section \ref{sec:construction} to
compute $d(q,A)$ and record a uniform convergence estimate for
Remark \ref{rem:varying-parameters}.  A Fourier estimate handles
all fields of order at least $49$; Cauchy--Davenport and Kneser's
theorem complete the small-field comparisons.

\subsection{Expected densities and uniform convergence}\label{sec:densities}

For a fixed nice pair $(q,A)$, formula \eqref{eq:finite-daisy}
from Section \ref{sec:construction} gives the expected density.
Its numerator product tends to $\beta_q$, while
$1-rq^{-r}\le(1-q^{-r})^r\le1$.  Hence
\[
 d(q,A)=\lim_{r\to\infty}p_r(q,A)
 =\beta_q\frac{|A|}{q-1}.
\]

For the variant in Remark \ref{rem:varying-parameters}, we record
an error bound that holds simultaneously for every nice pair.

\begin{lemma}[Uniform convergence of the expected densities]\label{lem:uniform-density}
For every $(q,A)\in\mathcal U$ and every $r\ge2$,
\begin{equation}\label{eq:uniform-density-error}
 0\le p_r(q,A)-d(q,A)\le(r+1)2^{-r}.
\end{equation}
In particular, the convergence is uniform over $\mathcal U$.
\end{lemma}

\begin{proof}
We denote the finite product and the nonzero-label correction factor by
\[
 T_r(q):=\prod_{k=1}^{r}(1-q^{-k}),
 \qquad c_r(q):=(1-q^{-r})^r.
\]
By \eqref{eq:finite-daisy},
\[
 p_r(q,A)=\frac{|A|}{q-1}\frac{T_r(q)}{c_r(q)},
 \qquad d(q,A)=\frac{|A|}{q-1}\beta_q.
\]
We estimate separately the tail of the infinite product and the
correction for excluding zero labels.  For the tail, the elementary
product inequality $1-\prod_j(1-u_j)\le\sum_j u_j$, with
$0\le u_j\le1$, gives
\[
 \begin{aligned}
 0\le T_r(q)-\beta_q
 &=T_r(q)\left(1-\prod_{k>r}(1-q^{-k})\right)\\
 &\le\sum_{k>r}q^{-k}
 =\frac{q^{-r}}{q-1}.
 \end{aligned}
\]
Next, $T_r(q)/c_r(q)$ is the probability of obtaining a basis from
independent uniform nonzero labels, and hence is at most $1$.
Consequently, Bernoulli's inequality gives
\[
 \begin{aligned}
 0\le\frac{T_r(q)}{c_r(q)}-T_r(q)
 &=\frac{T_r(q)}{c_r(q)}\bigl(1-c_r(q)\bigr)\\
 &\le1-(1-q^{-r})^r
 \le rq^{-r}.
 \end{aligned}
\]
Since $|A|/(q-1)\le1$ and $q\ge2$, combining the two estimates yields
\[
 0\le p_r(q,A)-d(q,A)
 \le rq^{-r}+\frac{q^{-r}}{q-1}
 \le(r+1)2^{-r}.
\]
The final bound is independent of both $q$ and $A$ and tends to zero.
\end{proof}

For $M_r$ in Remark \ref{rem:varying-parameters}, taking suprema
in \eqref{eq:uniform-density-error} gives
\begin{equation}\label{eq:uniform-supremum}
 \sup_{(q,A)\in\mathcal U}d(q,A)
 \le M_r
 \le\sup_{(q,A)\in\mathcal U}d(q,A)+(r+1)2^{-r}.
\end{equation}
This identifies the limit of $M_r$ with the fixed-rule supremum.
We now determine that supremum by bounding $d(q,A)$ over all nice
pairs in the remaining subsections.

\subsection{A Fourier bound over finite fields}

We work over an arbitrary finite field, including prime fields.
The first ingredient is a uniform estimate obtained by combining
the additive and multiplicative structures of the filter.

We fix a nontrivial additive character $\psi:\F_q\to\mathbb C$.
Thus $|\psi(x)|=1$, $\psi(x+y)=\psi(x)\psi(y)$, and
$\psi(-x)=\overline{\psi(x)}$.  The standard orthogonality relation
\cite[Section~4.1, in particular Lemma~4.5]{TaoVu} gives
\[
 \frac1q\sum_{t\in\F_q}\psi(tz)
 =\begin{cases}
 1,&z=0,\\
 0,&z\ne0.
 \end{cases}
\]
We use this identity to count solutions of the forbidden equation.

\begin{lemma}[A Fourier bound]\label{lem:fourier-bound}
Let $q$ be a prime power and let
$A\subseteq\F_q^\times$ be nonempty, sign-symmetric, and nice.  Then
\[
 |A|\le q^{3/4}.
\]
\end{lemma}

\begin{proof}
Let $m=|A|$, and let $N$ count the ordered six-tuples in $A^6$
satisfying
\begin{equation}\label{eq:six-variable}
 ab+cd=ef,\qquad a,b,c,d,e,f\in A,
\end{equation}
where repetitions are allowed.  By Theorem \ref{thm:criterion},
\eqref{eq:six-variable} has no solutions, so $N=0$.  We define
\[
 S(t):=\sum_{a,b\in A}\psi(tab).
\]
Applying the zero-detection identity to $ab+cd-ef$, summing over
$A^6$, and factoring the character sum gives
\[
 0=N=\frac1q\sum_{t\in\F_q}
 \sum_{a,b,c,d,e,f\in A}\psi\bigl(t(ab+cd-ef)\bigr)
 =\frac1q\sum_{t\in\F_q}S(t)^2S(-t).
\]
Since $S(0)=m^2$, separating $t=0$ yields
\[
 m^6=-\sum_{t\ne0}S(t)^2S(-t).
\]
Using $S(-t)=\overline{S(t)}$ and taking absolute values, we obtain
\begin{equation}\label{eq:fourier-main}
 m^6\le\sum_{t\ne0}|S(t)|^3
 \le\left(\max_{t\ne0}|S(t)|\right)
       \sum_{t\in\F_q}|S(t)|^2.
\end{equation}
We estimate the two factors separately.  First, for $t\ne0$,
Cauchy--Schwarz and orthogonality give
\[
 \begin{aligned}
 |S(t)|^2
 &\le m\sum_{a\in A}\left|\sum_{b\in A}\psi(tab)\right|^2
 \le m\sum_{a\in\F_q}\left|\sum_{b\in A}\psi(tab)\right|^2\\
 &=m\sum_{b,b'\in A}\sum_{a\in\F_q}
       \psi\bigl(ta(b-b')\bigr)=qm^2.
 \end{aligned}
\]
The inner sum is $q$ when $b=b'$ and zero otherwise.  Hence
$\max_{t\ne0}|S(t)|\le m\sqrt q$.
Second, expanding the square and applying orthogonality over $t$ gives
\[
 \begin{aligned}
 \sum_{t\in\F_q}|S(t)|^2
 &=\sum_{a,b,c,d\in A}\sum_{t\in\F_q}
       \psi\bigl(t(ab-cd)\bigr)\\
 &=q\,|\{(a,b,c,d)\in A^4:ab=cd\}|\le qm^3.
 \end{aligned}
\]
Indeed, once $a,b,c\in A$ are chosen, $c\ne0$ uniquely determines
$d=ab/c$, so there are at most $m^3$ such quadruples.
Substituting both bounds into \eqref{eq:fourier-main} gives
$m^6\le(m\sqrt q)(qm^3)=m^4q^{3/2}$, and hence $m\le q^{3/4}$.
\end{proof}

\subsection{Small fields and completion of the proof}\label{sec:small-fields}

We first record the comparison
\begin{equation}\label{eq:beta5-lower}
 \begin{aligned}
 \beta_5
 &=\frac45\prod_{k\ge2}(1-5^{-k})\\
 &>\frac45\left(1-\sum_{k\ge2}5^{-k}\right)
 =\frac{19}{25}.
 \end{aligned}
\end{equation}
For prime fields, we use Cauchy--Davenport; the remaining small
extension fields will be treated with Kneser's theorem.

\smallskip
\Needspace{8\baselineskip}
\noindent\textbf{Prime fields.}\par\nobreak
\begin{theorem}[Cauchy--Davenport]\label{thm:cauchy-davenport}
Let $p$ be prime and let $X,Y\subseteq\F_p$ be nonempty.  Then
\[
 |X+Y|\ge\min\{p,\,|X|+|Y|-1\}.
\]
\end{theorem}

A proof can be found in \cite[Theorem~5.4]{TaoVu}.

\begin{lemma}\label{lem:prime-bound}
Let $p$ be prime and let $A\subseteq\F_p^\times$ be nonempty, sign-symmetric,
and nice.  Then
\begin{equation}\label{eq:A-bound}
 |A|\le \left\lfloor\frac{p+1}{3}\right\rfloor.
\end{equation}
\end{lemma}

\begin{proof}
We choose $a_0\in A$ and replace $A$ by $a_0^{-1}A$.  This scaling preserves
cardinality, sign symmetry, and niceness, and it ensures that $1\in A$.
For $B=AA$ we then have $A\subseteq B$.  Theorem \ref{thm:criterion} says that $B+B$ and $B$ are disjoint.  The
Cauchy--Davenport theorem (Theorem \ref{thm:cauchy-davenport}), applied
with $X=Y=B$, gives
\[
 |B+B|\ge\min\{p,2|B|-1\}.
\]
The first alternative is impossible because $B$ is nonempty and disjoint
from $B+B$.  Therefore
\[
 p\ge |B+B|+|B|\ge3|B|-1,
\]
and hence
\[
 |A|\le|B|\le\left\lfloor\frac{p+1}{3}\right\rfloor.
\]
\end{proof}

\begin{proposition}\label{prop:prime-density}
For every prime $p$ and every $(p,A)\in\mathcal U$,
\[
 d(p,A)\le\frac{\beta_5}{2}.
\]
Equality holds precisely when $p=5$ and $A=\{\pm a\}$ for some
$a\in\F_5^\times$.
\end{proposition}

\begin{proof}
For $p=5$, every set $A=\{\pm a\}$ is nice: after factoring out
$a^2$, a sum of three elements of $AA=\{\pm a^2\}$ is one of
$\pm1,\pm3$ and is therefore nonzero in $\F_5$.  Lemma
\ref{lem:prime-bound} shows that no larger $A$ is possible.

We compare this choice with the remaining primes using \eqref{eq:beta5-lower}.
For $p=2$, the only nonempty filter is $A=\{1\}$, and
\[
 d(2,A)=\beta_2<\frac38<\frac{19}{50}<\frac{\beta_5}{2}.
\]
There is no nonempty nice filter in characteristic $3$, since any
$b\in AA$ satisfies $b+b+b=0$.  For $p=7$, Lemma
\ref{lem:prime-bound} gives $|A|/(p-1)\le1/3$.  For $p=11$, it gives
$|A|\le4$, and therefore
\[
 d(11,A)<\frac{10}{11}\frac4{10}=\frac4{11}.
\]
If $p\ge13$, we write $p=6k\pm1$.  The same lemma gives
\[
 \frac{|A|}{p-1}\le
 \begin{cases}
  1/3,&p=6k+1,\\[2mm]
  k/(3k-1)\le3/8,&p=6k-1,
 \end{cases}
\]
where the second case has $k\ge3$.  Since $3/8<19/50$, the limiting
densities in all these cases are strictly below $\beta_5/2$.

\end{proof}

\smallskip
\noindent\textbf{Small extension fields.}\par\nobreak
We use the following standard sumset theorem, stated in multiplicative
notation.

\begin{theorem}[Kneser]\label{thm:kneser}
Let $G$ be a finite abelian group, written multiplicatively, and let
$X,Y\subseteq G$ be nonempty.  Let $XY$ denote their product set and
$H$ its stabilizer under multiplication:
\[
 XY=\{xy:x\in X,\ y\in Y\},
 \qquad H=\{g\in G:gXY=XY\}.
\]
Then $H$ is a subgroup of $G$, and
\[
 \begin{aligned}
 |XY|&\ge |XH|+|YH|-|H|\\
     &\ge |X|+|Y|-|H|.
 \end{aligned}
\]
\end{theorem}

A proof can be found in \cite[Theorem~5.5]{TaoVu}, where the theorem
is stated in additive notation; see also DeVos \cite{DeVos} for a short
proof based on an intersection--union argument.  We apply it to the multiplicative
group $\F_q^\times$, not to the additive group of the field.

\begin{proof}[Completion of the proof of Theorem \ref{thm:finite-field-optimal}]
Let $q$ be any prime power.  Since
\[
 \beta_q
 =\prod_{k\ge1}(1-q^{-k})
 <1-q^{-1}=\frac{q-1}{q},
\]
Lemma \ref{lem:fourier-bound} implies
\begin{equation}\label{eq:large-q-bounds}
 d(q,A)<\frac{|A|}{q}\le q^{-1/4}.
\end{equation}
By \eqref{eq:beta5-lower},
\[
 \frac{\beta_5}{2}>\frac{19}{50}>\frac1{\sqrt7}.
\]
Consequently \eqref{eq:large-q-bounds} proves the strict inequality for
every $q\ge49$.

It remains to treat $q<49$.  For prime $q$, Proposition
\ref{prop:prime-density} gives
$d(q,A)\le\beta_5/2$, with equality precisely when $q=5$ and
$A=\{\pm a\}$ for some $a\in\F_5^\times$.
The remaining proper prime powers below $49$ are
\[
 4,8,9,16,25,27,32.
\]
There is no nice filter in characteristic $3$: if $b\in B:=AA$, then
$b+b=-b\in B$ because $B=-B$, contradicting
$(B+B)\cap B=\varnothing$.  This excludes $q=9$ and $q=27$.

For the remaining fields, first observe that every nonempty additively
sum-free set $B$ satisfies
\begin{equation}\label{eq:B-half}
 |B|\le\left\lfloor\frac q2\right\rfloor.
\end{equation}
Indeed, for any $b\in B$, the sets $B$ and $B+b$ are disjoint.  We normalize
$A$ by a nonzero scalar so that $1\in A$; this preserves niceness and
ensures that $A\subseteq B=AA$.

The small-field bounds needed below are
\begin{equation}\label{eq:small-field-bounds}
\begin{array}{c|ccccc}
 q&4&8&16&25&32\\ \hline
 |A|\text{ is at most}&1&2&6&9&8.
\end{array}
\end{equation}
We prove these bounds using Kneser's theorem, under the normalization
$1\in A$ and the bound \eqref{eq:B-half}.

\smallskip
\noindent\textbf{The small-field bounds.}
Let $H$ be the stabilizer of $B=AA$ in the multiplicative group
$\F_q^\times$, and let $h=|H|$.  Kneser's theorem
(Theorem \ref{thm:kneser}), applied with $X=Y=A$, gives
\begin{equation}\label{eq:mult-kneser}
 |B|\ge2|A|-h.
\end{equation}
Moreover, $h\mid(q-1)$ and $h\mid|B|$, since $B$ is a union of
$H$-cosets.  In particular, $h\le|B|\le\lfloor q/2\rfloor$.
For $q=4,8,32$, the group order $q-1$ is prime, so these restrictions
force $h=1$.  Hence
\[
 2|A|-1\le\left\lfloor\frac q2\right\rfloor,
\]
which gives $|A|\le1,2,8$, respectively.
For $q=16$, the conditions $h\mid15$ and $h\le8$ give
$h\in\{1,3,5\}$.  Thus
\[
 2|A|\le|B|+h\le8+5=13,
\]
and $|A|\le6$.

For $q=25$, suppose that $|A|\ge10$.  Then
\[
 20-h\le|B|\le12,
\]
so $8\le h\le12$.  As $h\mid24$, only $h=8,12$ are possible.
If $h=8$, then $|B|\ge12$ and $8\mid|B|$, forcing $|B|\ge16$,
a contradiction.  Hence $h=|B|=12$.  Since $1\in B$, we have $B=H$,
the unique subgroup of order $12$ in the cyclic group
$\F_{25}^\times$, namely the subgroup of squares.  The element
$2\in\F_5^\times$ satisfies $2^{12}=1$, so it belongs to this subgroup.
Consequently $1,2\in B$ and $1+1=2$, contrary to the sum-free property.
Thus $|A|\le9$, proving \eqref{eq:small-field-bounds}.

\smallskip
\noindent\textbf{Completing the comparison.}
For these five fields, \eqref{eq:small-field-bounds} gives
\[
 \frac{|A|}{q}\le
 \frac14,\frac14,\frac38,\frac9{25},\frac14,
\]
respectively.  Hence \eqref{eq:large-q-bounds} (without using the
$q^{-1/4}$ estimate) and \eqref{eq:beta5-lower} show that $d(q,A)$
is strictly below $\beta_5/2$.  Thus
\[
 \max_{(q,A)\in\mathcal U}d(q,A)=\frac{\beta_5}{2},
\]
with the equality cases given by the prime-field argument.
This completes the proof of Theorem \ref{thm:finite-field-optimal}.

\end{proof}

\section{Applications with improved constants}\label{sec:applications}

\subsection{An induced \texorpdfstring{$C_6$}{C6}-free graph in one hypercube layer}

For $r\ge2$, fix a nonzero vector $v_0\in\F_q^r$ and choose
$v_1,\dots,v_n$ independently and uniformly from
$\F_q^r\setminus\{0\}$.  We retain an upper $r$-set $S$ when
$\Delta(S)\in A$, and a lower $(r-1)$-set $T$ when
$[v_0,v_i:i\in T]\in A$.  Let $G_r(q,A)$ be the subgraph of
$L_r(n)$ induced by the retained sets.

The number of ordered completions of $v_0$ to a basis is
\[
 \prod_{j=1}^{r-1}(q^r-q^j)
 =\frac{|\mathrm{GL}_r(q)|}{q^r-1}.
\]
Since at least one random column can be scaled, each nonzero
determinant value occurs $|\mathrm{SL}_r(q)|/(q^r-1)$ times.
Dividing by the $(q^r-1)^{r-1}$ possible samples gives
\begin{equation}\label{eq:fixed-det-prob}
 \Prob\bigl([v_0,v_1,\dots,v_{r-1}]\in A\bigr)
 =\frac{|A|\,|\mathrm{SL}_r(q)|}{(q^r-1)^r}.
\end{equation}
Conditional on a lower vertex $T$ being retained, its columns are
linearly independent.  Thus $x\mapsto[v_i:i\in T,x]$ is a nonzero
linear functional, with $q^{r-1}$ nonzero preimages of each prescribed
nonzero value.  The conditional probability of retaining
$T\cup\{j\}$ is therefore $|A|q^{r-1}/(q^r-1)$.  By linearity of
expectation,
\begin{equation}\label{eq:layer-limit}
 \begin{aligned}
 \frac{\E e(G_r(q,A))}{e(L_r(n))}
 &=\frac{|A|^2q^{r-1}|\mathrm{SL}_r(q)|}{(q^r-1)^{r+1}},\\
 \lim_{r\to\infty}\frac{\E e(G_r(q,A))}{e(L_r(n))}
 &=\beta_q\frac{|A|^2}{q(q-1)}.
 \end{aligned}
\end{equation}
Here $(q,A)$ is fixed as $r\to\infty$, while $n\ge r$ may vary
arbitrarily with $r$.

\begin{proof}[Proof of Corollary \ref{cor:layer-main}]
We specialize the construction above to $q=5$, $A=\{\pm1\}$, and
write $G_r=G_r(5,\{\pm1\})$.
Suppose $G_r$ contains a $6$-cycle.  Every such cycle in a Boolean
edge-layer has a common $(r-2)$-set $I$ and three further indices
$a,b,c$ such that the lower vertices are
\[
I\cup\{a\},\quad I\cup\{b\},\quad I\cup\{c\},
\]
and the upper vertices are
\[
I\cup\{a,b\},\quad I\cup\{a,c\},\quad I\cup\{b,c\}.
\]
We order the common $I$-columns first and regard $v_0,v_a,v_b,v_c$
as the four remaining columns.  In the common-stem brackets, we
abbreviate $[v_i v_j]$ to $[ij]$.  Membership of the six vertices implies
\[
[0a],[0b],[0c],[ab],[ac],[bc]\in\{\pm1\}.
\]
The Grassmann--Pl\"ucker relation
\[
[0a][bc]-[0b][ac]+[0c][ab]=0
\]
is impossible in $\F_5$ for the same reason as in the proof of Theorem \ref{thm:daisy-main}.  Hence no such $C_6$ exists.

By \eqref{eq:layer-limit} and \eqref{eq:SL},
\[
 \frac{\E e(G_r)}{e(L_r(n))}
 =\frac{4\cdot5^{r-1}|\mathrm{SL}_r(5)|}{(5^r-1)^{r+1}}
 =\frac15\frac{\prod_{k=1}^{r-1}(1-5^{-k})}{(1-5^{-r})^r}
 >\frac{\beta_5}{5}.
\]
Thus some outcome is induced $C_6$-free and has more than
$(\beta_5/5)e(L_r(n))$ edges.
\end{proof}

\Needspace{8\baselineskip}
\smallskip
\noindent\textbf{Optimality within the two-level construction.}
The limiting value $\beta_5/5$ is also optimal among nice filters.

\begin{proposition}\label{prop:layer-optimal}
Let $q$ be a prime power and let $A\subseteq\F_q^\times$ be nonempty,
sign-symmetric, and nice.  The two-level construction
above satisfies
\[
 \lim_{r\to\infty}\frac{\E e(G_r(q,A))}{e(L_r(n))}
 =\beta_q\frac{|A|^2}{q(q-1)}
 \le\frac{\beta_5}{5}.
\]
Equality holds precisely when $q=5$ and $A=\{a,-a\}$ for some
$a\in\F_5^\times$.
\end{proposition}

\begin{proof}
The density formula is \eqref{eq:layer-limit}.  We reuse the bounds on
$|A|$ proved in Section \ref{sec:optimality}, and compare the resulting
densities with
\[
 \frac{\beta_5}{5}>\frac{19}{125},
\]
which follows from \eqref{eq:beta5-lower}.

First suppose that $q=p$ is prime.  For $p=5$, Lemma
\ref{lem:prime-bound} gives $|A|\le2$, and every sign-symmetric pair
$\{a,-a\}$ is nice and gives equality.  For $p=2$, necessarily
$A=\{1\}$, and
\[
 \frac{\beta_2}{2}
 <\frac12\prod_{j=1}^{5}(1-2^{-j})
 =\frac{9765}{65536}<\frac{19}{125}.
\]
There is no nonempty nice filter in characteristic $3$.
For $p=7$ and $p=11$, respectively, Lemma \ref{lem:prime-bound} gives
\[
 \beta_7\frac{|A|^2}{7\cdot6}<\frac{2}{21},
 \qquad
 \beta_{11}\frac{|A|^2}{11\cdot10}<\frac{8}{55},
\]
both less than $19/125$.  For $p\ge13$, we write $p=6k\pm1$.
If $p=6k+1$, both $|A|/(p-1)$ and $|A|/p$ are at most $1/3$.
If $p=6k-1$, then $k\ge3$, and these factors are at most $3/8$
and $6/17$, respectively.  Thus in either case
\[
 \beta_p\frac{|A|^2}{p(p-1)}
 <\frac9{68}<\frac{19}{125}.
\]

Now suppose that $q$ is a proper prime power.  Since
$\beta_q<(q-1)/q$, the Fourier bound (Lemma \ref{lem:fourier-bound})
gives, for $q\ge49$,
\[
 \beta_q\frac{|A|^2}{q(q-1)}
 <\frac{|A|^2}{q^2}
 \le q^{-1/2}
 \le\frac17<\frac{19}{125}.
\]
Among the smaller proper prime powers, characteristic $3$ excludes
$q=9,27$.  For $q=4,8,16,25,32$, the bounds in
\eqref{eq:small-field-bounds} give $|A|/q\le3/8$, so
\[
 \beta_q\frac{|A|^2}{q(q-1)}
 <\left(\frac{|A|}{q}\right)^2
 \le\frac9{64}<\frac{19}{125}.
\]
All cases other than the stated $\F_5$ pairs therefore give strict
inequality.
\end{proof}

This proposition concerns only the limiting expected density of this
two-level construction.  It does not assert that $\beta_5/5$ is an upper
bound for arbitrary induced $C_6$-free subgraphs of a Boolean edge-layer.

\subsection{A four-dimensional local-density consequence}

We apply Theorem \ref{thm:daisy-main} to the first nontrivial case of the
local-density question of Johnson and Talbot \cite{JT}, following the
construction of Ellis, Ivan and Leader \cite[Theorem~8]{EIL}, which
selects every third Boolean layer.

\begin{proof}[Proof of Corollary \ref{cor:q4}]
We identify the vertex set of $Q_n$ with the power set
$\mathcal P([n])$.  For each $3\le r\le n$ with
$r\equiv0\pmod3$, we choose a $D_r$-free family
\[
\calF_r\subseteq\binom{[n]}r
\]
with density at least $\betafive/2$, as supplied directly by Theorem
\ref{thm:daisy-main}.  We set $\calF_0=\varnothing$ and define
\[
\mathcal A=\bigcup_{\substack{0\le r\le n\\r\equiv0\ (3)}}\calF_r,
\]

Let $\omega=e^{2\pi i/3}$.  The roots-of-unity filter gives, for
$j\in\{0,1,2\}$,
\[
 \begin{aligned}
 M_j(n)&:=\sum_{\substack{0\le r\le n\\r\equiv j\ (3)}}\binom nr\\
 &=\frac13\left(2^n+\omega^{-j}(1+\omega)^n
       +\omega^{-2j}(1+\omega^2)^n\right)
 =\frac{2^n}{3}+O(1),
 \end{aligned}
\]
since $|1+\omega|=|1+\omega^2|=1$.  Omitting the zeroth layer
changes the normalized size by only $O(2^{-n})$, so
\[
 \frac{|\mathcal A|}{2^n}
 \ge\frac{\beta_5}{2}\left(\frac13-O(2^{-n})\right)
 =\frac{\beta_5}{6}-o(1).
\]

Let $V$ be the vertex set of any $4$-dimensional subcube.  It meets five consecutive Boolean layers, say $r_0,r_0+1,\dots,r_0+4$, with multiplicities $1,4,6,4,1$.

If $r_0\equiv0\pmod3$, then only the layers $r_0$ and $r_0+3$ can contribute to $\mathcal A$, giving at most $1+4=5$ points.  If $r_0\equiv2\pmod3$, only $r_0+1$ and $r_0+4$ contribute, again giving at most $4+1=5$ points.  Finally, if $r_0\equiv1\pmod3$, the only selected layer is the middle layer $r_0+2$, where the six vertices of $V$ form an $(r_0+2)$-daisy.  Since $\calF_{r_0+2}$ is daisy-free, at most five of those six vertices belong to $\mathcal A$.  This proves the claim.
\end{proof}

\section{Concluding remarks}\label{sec:conclusion}

\begin{remark}[Scope of the criterion]\label{rem:criterion-scope}
Theorem \ref{thm:criterion} characterizes rules that guarantee
daisy-freeness for every labelling; its algebraic condition need not
hold for a particular daisy-free outcome.  For example, with $r=2$,
$q=3$, and $A=\{\pm1\}$, the condition fails because
$1+1=-1\in AA$.  If all labels belong to $\{e_1,e_2\}$, however,
the retained pairs form a complete bipartite graph and hence contain
no $D_2=K_4$.  Similarly, Theorem \ref{thm:finite-field-optimal}
optimizes the limiting expected density within the stated framework;
it gives no upper bound on $\pi(D_r)$.
\end{remark}

\subsection{Nonuniform sampling}

One possible direction beyond Theorem \ref{thm:finite-field-optimal}
is to relax the assumption that the vectors are chosen uniformly,
while retaining independence and a common distribution.  For any
nice determinant filter, the daisy-free conclusion
holds for every labelling, by Theorem \ref{thm:criterion}.  Uniformity
is used to compute and optimize the expected density, not to exclude
a daisy.

For a prime power $q$, a nice filter
$A\subseteq\F_q^\times$, and a probability distribution $\mu$ on
$\F_q^r\setminus\{0\}$, let $P_r(q,A;\mu)$ denote the probability
that $r$ independent samples from $\mu$ have an accepted determinant:
\[
 P_r(q,A;\mu):=
 \Prob_{v_1,\ldots,v_r\,\mathrm{iid}\sim\mu}
 \bigl([v_1,\ldots,v_r]\in A\bigr).
\]
Assigning labels independently with this common distribution gives a
daisy-free family of expected density $P_r(q,A;\mu)$.  We denote the
uniform distribution on $\F_q^r\setminus\{0\}$ by $u_{q,r}$.
For this distribution, equation \eqref{eq:random-det-prob} gives
\[
 P_r(q,A;u_{q,r})=\frac{|A|}{q-1}\,
 \frac{\prod_{k=1}^r(1-q^{-k})}{(1-q^{-r})^r}
 \longrightarrow d(q,A).
\]

A small-dimensional example shows that uniformity need not maximize
this probability.  We take $q=5$, $A=\{\pm1\}$, and $r=2$, and let
$\mu$ assign mass $1/3$ to each of $e_1,e_2,e_1+e_2$ and zero to all
other nonzero vectors.  Any two distinct vectors in this support have
determinant in $\{\pm1\}$, whereas a repeated vector gives determinant
zero.  Thus
\[
 P_2(5,\{\pm1\};\mu)=\frac23
 >\frac5{12}=P_2(5,\{\pm1\};u_{5,2}).
\]
In particular, the analogue of Alon's fixed-dimensional
uniform-optimality statement for the binary basis construction
\cite[Corollary~2.4]{Alon} fails for this determinant filter.

The asymptotic question is whether relaxing uniformity can improve
the limiting lower bound for the Tur\'an density of daisies, rather than merely improve
the density in a fixed dimension.  More precisely, can one choose
prime powers $q_r$, nice filters
$A_r\subseteq\F_{q_r}^\times$, and distributions $\mu_r$ on
$\F_{q_r}^r\setminus\{0\}$ such that
\begin{equation}\label{eq:nonuniform-question}
 \liminf_{r\to\infty}P_r(q_r,A_r;\mu_r)>\frac{\beta_5}{2}?
\end{equation}
The field and filter may in particular be fixed; the labels remain
independent and identically distributed at each rank.  The example
above motivates this question but does not establish an asymptotic
improvement.

\subsection{Further consequences of the Ellis--Ivan--Leader framework}

Ellis, Ivan and Leader \cite{EIL} use a positive-density ordinary
daisy-free family as an input to several constructions.  Wherever an
argument uses only that property and its density, Theorem
\ref{thm:daisy-main} can replace their binary family, changing the input
constant from $\beta_2$ to $\beta_5/2$ while leaving the rest of the
argument intact.  Corollary \ref{cor:q4} gives one such application.
We record two further consequences here without repeating the proofs.

First, in the general local-density construction in
\cite[Theorem~8]{EIL}, we fix $d\ge5$ and let
$m=\lceil C\sqrt d\rceil$, with the same sufficiently large absolute
constant $C$ as in their proof.  We insert our daisy-free family in each
Boolean layer whose index is a multiple of $m$.  Their argument that
every $d$-dimensional subcube contains fewer than
$\binom d{\lfloor d/2\rfloor}$ selected vertices is unchanged, while
the resulting family $\mathcal A$ satisfies
\[
 \frac{|\mathcal A|}{2^n}
 \ge \frac{\beta_5}{2m}-o(1)
 \qquad(n\to\infty).
\]
Thus the numerator $\beta_2$ in their density estimate is replaced by
$\beta_5/2$; the order $d^{-1/2}$ and the layer spacing are not improved.

Second, let $B_d$ be the Boolean lattice of all subsets of a $d$-element
set, ordered by inclusion.  In their counterexample to the layer-union
conjecture of Bukh and of Griggs and Lu, Ellis, Ivan and Leader
\cite{EIL} take $d+1$ consecutive middle layers and replace a middle
layer by an ordinary daisy-free family.  For each fixed $d\ge4$,
their forbidden-subposet argument applies to our family and gives
a $B_d$-free family $\mathcal G\subseteq\mathcal P([n])$ with
\[
 |\mathcal G|
 \ge \left(d+\frac{\beta_5}{2}-o(1)\right)
       \binom n{\lfloor n/2\rfloor},
\]
in place of the coefficient $d+\beta_2$ supplied by the binary family.
Here containment is as a subposet, not necessarily an induced subposet.

These are consequences of their existing arguments, with a stronger
bound for the density of an ordinary daisy-free family.  This
substitution principle should not be read as applying to every result
of \cite{EIL}: their generalized-daisy constructions, used for example
in their exponential bounds for hypercube hitting sets, require
different forbidden configurations and are not supplied by our
ordinary-daisy theorem.

\subsection{What the construction gives for \texorpdfstring{$C_{10}$}{C10}}

For a graph $H$, let $\ex(Q_n,H)$ denote the maximum number of edges
in an $H$-free subgraph of $Q_n$.  We denote its hypercube Tur\'an
density by
\[
 \pi_\square(H):=\lim_{n\to\infty}
 \frac{\ex(Q_n,H)}{e(Q_n)}.
\]

The induced-$C_6$-free layer in Corollary \ref{cor:layer-main} can also
be substituted into the construction of Grebennikov and Marciano
\cite{GMv4}.  Their argument uses the Axenovich--Martin two-colouring
of a Boolean edge-layer with no monochromatic induced $C_{10}$
\cite{AM}.  Keeping a larger colour class retains at least half the
edges of our layer graph and makes it $C_{10}$-free: every $C_{10}$
in an induced $C_6$-free subgraph of the layer must be induced in the
layer, since a chord would create a $C_4$ or a $C_6$.
Taking alternate edge-layers incurs another asymptotic factor of
$1/2$ and prevents cycles from joining different retained layers.
Consequently, the same argument gives
\[
 \begin{aligned}
 \pi_\square(C_{10})
 &\ge \frac14\cdot\frac{\beta_5}{5}
 =\frac{\beta_5}{20}\\
 &=0.03801663979\ldots,
 \end{aligned}
\]
in place of the constant
\[
 \frac{\beta_2}{8}=0.0360985\ldots
\]
obtained by that route with their binary layer construction.  Here we
refer specifically to the revised arXiv version \cite{GMv4}, rather
than the earlier journal version, whose stated constant is
\[
 \frac{\beta_2}{12}=0.0240656\ldots.
\]

This comparison concerns the Grebennikov--Marciano framework, not an
improvement of Marko Peji\'c's result.  Peji\'c
\cite[Theorem~1.2]{Pejic} proved the stronger bound
\[
 \pi_\square(C_{10})\ge\frac18.
\]
His proof uses a different edge-retention rule over $\F_3$, together
with repeated independent rounds and deletion of previously selected
vertices, followed by the selection of alternate edge-layers. 
Since
\[
 \frac{\beta_5}{20}=0.03801663979\ldots<\frac18,
\]
we include our bound only to illustrate this further application
of the determinant construction, rather than as a main theorem or
corollary.

\clearpage

\section*{Disclosure of AI assistance}

The results in this paper were developed through extended conversations
with GPT-5.6 Sol and GPT-6 Astra, primarily concerning the results and
constructions of Ellis, Ivan and Leader \cite{EIL}.
More precisely, the AI proposed the determinant-value construction over
$\F_5$.  After understanding this construction, the author formulated
the general optimization problem .
This problem was subsequently solved through further discussions
between the author and the AI.  The author independently reconstructed
the proofs, revised the exposition, and takes full responsibility for
the mathematical content.

\end{document}